\documentclass[11pt]{article}
\usepackage[centertags]{amsmath}
\usepackage{amsfonts}
\usepackage{amssymb}
\usepackage{amsthm}
\usepackage{newlfont}
\usepackage{graphicx}
\usepackage{multirow}
\usepackage{rotating}
\usepackage{epstopdf}

\newtheorem{theor}{Theorem}
\newtheorem{definition}{Definition}

\newtheorem{example}{Example}

\begin{document}
\thispagestyle{empty}
%
%
\title{A novel technique to solve the fuzzy system of equations}

\author{ N. Mikaeilvand$^{a,}$\footnote{ E-mail:
Mikaeilvand@iauardabil.ac.ir} , Z. Noeiaghdam$^{b,}$\footnote{
E-mail: Z$_{-}$Noeiaghdam@shahed.ac.ir; zahra.noie@yahoo.com } , S.
Noeiaghdam$^{c,d,}$\footnote{ Corresponding author, E-mail:
s.noeiaghdam.sci@iauctb.ac.ir; samadnoeiaghdam@gmail.com} , Juan J.
Nieto$^{e,}$\footnote{ E-mail: juanjose.nieto.roig@usc.es} }
\date{}
\maketitle \vspace{-9mm}
\begin{center}
\scriptsize{$^a$Department of Mathematics, Ardabil Branch, Islamic
Azad University, Ardabil, Iran.\\
$^b$Department of Mathematics and Computer Science, Shahed University, Tehran, Iran. \\
$^c$South Ural State University, Lenin prospect 76, Chelyabinsk, 454080, Russian Federation.\\
$^d$Baikal School of BRICS, Irkutsk National Research Technical
University, Irkutsk, Russian Federation.\\
$^e$Instituto de Matem\'aticas, Departamento de Estat\'istica,
An\'alise Matem\'atica e Optimizaci\'on, Universidade de Santiago de
Compostela 15782, Santiago de Compostela, Spain. }
\end{center}
\date{}

%
%
%
\maketitle

\begin{abstract}
The aim of this research is to apply a novel technique based on the
embedding method to solve the  $n\times n$ fuzzy system of linear
equations (FSLEs). By using this method, the strong fuzzy number
solutions of FSLEs can be obtained by transforming the
 $n \times n$ FSLE to the crisp system.
In this paper, Ezzati's method to solve the FSLEs is modified and
improved. Several theorems are proved to show the number of
operations for presented method are less than the methods of
Friedman and Ezzati. In order to show the advantages of scheme, two
applicable algorithms are presented and several examples are solved
by applying them. Also, some graphs of obtained results are
demonstrated which show the solutions are in the fuzzy form.

\vspace{.5cm}{\it Keywords:} Fuzzy linear system, Fuzzy number,
Fuzzy number vector, Embedding method.
\end{abstract}
\section{Introduction}
The FSLEs have many applications in different fields of science and
engineering such as heat transport, fluid flow, electromagnetism and
so on. In the recent decades, solving and studying the FSLEs have
been appeared in many researches. In 1998, Friedman et al. \cite{k1}
presented a model to solve the FSLEs and many mathematicians
improved and developed this method to find the solution of FSLEs. In
last years, Friedman et al. \cite{k1,k2,mm}, Abbasbany et al.
\cite{5,6,ba}, Allahviranloo et al. \cite{a5,a,ahkk,a6,aa6} and
others \cite{ma,e,FF,z1,z3,z2,kw,bz} considered the $n \times n$
FSLEs. Also, many authors applied the numerical methods to find the
approximate solution of the FSLEs \cite{5,6,karj,nag,SNE}.
Furthermore, the CESTAC method \cite{1,3,2,manham} based on the
stochastic arithmetic has been applied to find the optimal
iteration, optimal approximations and the optimal error of numerical
methods to solve the FSLEs \cite{4,7}.

In this research, the Ezzati's method \cite{e} to solve the $n
\times n$ FSLEs is improved and a novel method is presented. The aim
of this paper is to apply the embedding method and substituting the
 $n\times n$ FSLEs by two $n\times n$ crisp systems.
Several theorems and lemmas are proved that show the number of
operations in new method is lower than the ezzati's method and the
solutions of the FSLEs can be obtained by the fast and safe process.
In order to show the abilities of method, two algorithms are
presented and are applied to solve the examples.

Organization of this paper is in the following form: Section 2
contains several definitions and theorems of fuzzy arithmetic and
the FSLEs. Section 3 introduces the new idea to solve the FSLEs. In
this section, Ezzati's method is modified and improved. Also,
several theorems are proved to show the presented method is better
than the methods of Friedman and Ezzati. Furthermore, we will prove
the number of operations in presented method is less than the
mentioned methods. In Section 4, two applicable and efficient
algorithms are presented. Also, several examples are solved by using
these algorithms. Graphs of solutions are demonstrated to show the
results are in the fuzzy form. Section 5 is the conclusion.

\section{Preliminaries}
Several definitions and details of fuzzy arithmetic are presented
\cite{k1,mm}. Also, the methods of Friedman \cite{k1} and Ezzati
\cite{e} to solve the FSLEs are considered. Ezzati's method had some
problems in proving the theorems that we modify and represent them.
\begin{definition}\label{d2}
 \cite{k1,mm} Let $\widetilde{p} = (\underline{p}(z),\overline{p}(z)), \ 0
\leq z \leq 1 $ be the arbitrary fuzzy number then the following
criteria should be satisfied
\begin{itemize}
\item [(i)]$\underline{p}(z)$ is a bounded monotonic increasing
left continuous function.
 \item[(ii)]$\overline{p}(z)$ is a bounded monotonic decreasing
left continuous function.
 \item[(iii)]$\underline{p}(z)\leq
\overline{p}(z),\ 0 \leq z \leq 1. $
\end{itemize}
\end{definition}
The set of all fuzzy numbers is denoted by $\mathbf{E^{1}}$. The
crisp number $k$ is called the singleton when $\overline{p}(z)
=\underline{p}(z)=k ,~0 \leq z \leq 1$.

Let $\widetilde{p}=(\underline{p}(z),\overline{p}(z)),~
\widetilde{q}=(\underline{q}(z),\overline{q}(z))$  be the arbitrary
fuzzy functions and $k$ be the scalar value. The operations between
two fuzzy functions are defined as follows

$$(\underline{p+q})(z)=\underline{p}(z)+\underline{q}(z),~~~
(\overline{p+q})(z)=\overline{p}(z)+\overline{q}(z),$$

$$
(\underline{p-q})(z)=\underline{p}(z)-\overline{q}(z),~~~
(\overline{p-q})(z)=\overline{p}(z)-\underline{q}(z),$$

$$ \widetilde{k p} =\left\{\begin{array}{cc}
\displaystyle  (k\underline{p}(z),k\overline{p}(z)), & k\geq 0 , \\
  \\
\displaystyle  (k\overline{p}(z),k\underline{p}(z)), & k<0. \\
\end{array}
\right.
$$
Also $\widetilde{p}=\widetilde{q}$ if and only if
$\underline{p}(z)=\underline{q}(z) $ and
$\overline{p}(z)=\overline{q}(z)$.

\textbf{Remark 1:} The triangular fuzzy number
$\tilde{p}=(c,\mu,\rho)$ is defined as follows
\begin{equation}\label{scalar}
 \widetilde{p}(x)=\left\{\begin{array}{cc}
\displaystyle  \frac{x-c+\mu}{\mu}, &c-\mu\leq x\leq c, \\
  \\
\displaystyle  \frac{c+\rho-x}{\rho}, & c \leq x\leq c+\rho, \\
  \\
\displaystyle  0, &o.w,\\
\end{array}
\right.
\end{equation}
where~$\mu, \rho > 0$.~It is clear
that~$\underline{p}(z)=c-(1-z)\mu,~\overline{p}(z)=c+(1-z)\rho$ and
$\overline{p}-\underline{p}=(\mu+\rho)(1-z)$. The set of all
triangular fuzzy numbers is denoted by~$\mathbf{TE^{1}}.$

\begin{definition}\label{d4}
~\cite{k1} Let $(\tilde{v}_1,\tilde{v}_2,\ldots, \tilde{v}_j, \ldots
,\tilde{v}_n)^{T},
\tilde{v}_j=(\underline{v}_j(z),\overline{v}_j(z)),~1\leq j\leq
n,~0\leq z\leq 1$  be the fuzzy number vector which is called the
solution of FSLEs if and only if
\begin{equation}\label{3}
        \begin{array}{l}
   \displaystyle     \underline{\sum_{j=1}^na_{ij}v_j}(z)=\sum_{j=1}^n\underline{a_{ij}v_j}(z)=\underline{b_i}(z), \\
        \\
\displaystyle
\overline{\sum_{j=1}^na_{ij}v_j}(z)=\sum_{j=1}^n\overline{a_{ij}v_j}(z)=\overline{b_i}(z),~~~i=1,2,\cdots,n.
                \end{array}
      \end{equation}
\end{definition}

Finally, the methods of Friedman et al. \cite{k1} and Ezzati
\cite{e} to solve the FSLEs are reminded.

\subsection{Friedman's method}
Friedman et al. \cite{k1} presented the FSLEs as
\begin{equation}\label{sx=y}
Sv(z)=w(z)
\end{equation}
where
\begin{equation}\label{7}
    \begin{array}{l}
      v(z)=\Big(\underline{v}_1(z),\underline{v}_2(z),\cdots,\underline{v}_n(z),-\overline{v}_1(z),
      -\overline{v}_2(z),\cdots,-\overline{v}_n(z)\Big)^T,\\

      \\
      w(z)=\left(\underline{w}_1(z),\underline{w}_2(z),
      \cdots,\underline{w}_n(z),-\overline{w}_1(z),-\overline{w}_2(z),\cdots,-\overline{w}_n(z)\right)^T,
    \end{array}
\end{equation}
and the elements of $S = (s_{ij}),~1 \leq i, j \leq 2n$ are obtained
 based on the following conditions
\begin{equation}\label{6}
    \begin{array}{l}
      a_{ij} \geq 0 \Rightarrow s_{ij}=s_{i+nj+n}=a_{ij},\\
      \\
      a_{ij}<0 \Rightarrow s_{ij+n}=s_{i+nj}=-a_{ij}.
          \end{array}
\end{equation}
We note that for values $s_{ij}$ which are determined by neglecting
the criterion $(\ref{6})$ we have $s_{ij}=0$.  The matrix $S$ for
$s_{ij} \geq 0, 1\leq i,j \leq 2n$ can be formed as follows
\[S=\left(
      \begin{array}{cc}
        B & C \\
        C & B \\
      \end{array}
    \right),
\]
where $B$ constructs by the positive entries of $A$, $C$ constructs
by the absolute values of the negative entries of $A$ and $A = B -
C$.

For nonsingular matrix $S$ we have $v(z)=S^{-1}w(z)$. But probably
 the obtained solution has not the proper fuzzy number vector.
Therefore, the solution of the FSLEs can be defined in the following
form
\begin{definition}\label{d8}
\cite{k1} Let Eq. $(\ref{sx=y})$ has the unique solution in the form
 $\tilde{v}(z)=\{(\underline{v}_i(z),-\overline{v}_i(z)),~1\leq i\leq
n\}$. We define the fuzzy number vector
$\tilde{P}=\{(\underline{p}_i(z),\overline{p}_i(z)),~1\leq i\leq
n\}$ as
\[\underline{p}_i(z)=\min\{\underline{v}_i(z),\overline{v}_i(z),\underline{v}_i(1),\overline{v}_i(1)\},\]
\[\overline{p}_i(z)=\max\{\underline{v}_i(z),\overline{v}_i(z),\overline{v}_i(1),\underline{v}_i(1)\},\]
which is called the fuzzy solution of Eq. $(\ref{sx=y})$.
\end{definition}

If $\underline{p}_i(z)=\underline{v}_i(z)$ and
$~\overline{p}_i(z)=\overline{v}_i(z),~1\leq i\leq n$ then
$\tilde{P}$ is called a strong fuzzy solution. Otherwise,
$\tilde{P}$ is called a weak fuzzy solution which it is not the
solution of FSLE and it is not always fuzzy number vector. Recently,
Allahviranloo et al. \cite{a} showed that the weak solution of a
FSLE is not always a fuzzy number vector and it is the main fault of
Friedman's method.


\subsection{Ezzati's method \cite{e}}
Consider the following FSLEs
\begin{equation}\label{1}
 \left\{\begin{array}{ll}
a_{11}(\underline{v}_1(z)+\overline{v}_1(z))+\ldots+a_{1n}(\underline{v}_n(z)+\overline{v}_n(z))
&=\underline{w}_1(z)+\overline{w}_1(z), \\
\\
a_{21}(\underline{v}_1(z)+\overline{v}_1(z))+\ldots+a_{2n}(\underline{v}_n(z)+\overline{v}_n(z))&=
\underline{w}_2(z)+\overline{w}_2(z), \\
~~~~~~~~~~~~~~~~~~~~~~~~~~~~~~~~~\vdots & ~~~~~~~~~~~\vdots\\

a_{n1}(\underline{v}_1(z)+\overline{v}_1(z))+\ldots+a_{nn}(\underline{v}_n(z)+\overline{v}_n(z))
&=\underline{w}_n(z)+\overline{w}_n(z),\\
\end{array}
\right.
\end{equation}
where the solution of system (\ref{1}) is in the following form
\[\mathbf{g(z)}=\left(%
\begin{array}{c}
   g_{1}(z)\\
   g_{2}(z)\\
   \vdots\\
   g_{n}(z)\\
\end{array}%
\right)=
\mathbf{\underline{v}(z)+\overline{v}(z)}=\left(%
\begin{array}{cccc}
  \underline{v}_1(z)+\overline{v}_1(z)\\
   \underline{v}_2(z)+\overline{v}_2(z)\\
   \vdots\\
   \underline{v}_n(z)+\overline{v}_n(z)\\
\end{array}%
\right).\]

Since $(B+C)\underline{v}(z)=\underline{w}(z)+Cg(z)$ and
$(B+C)\overline{v}(z)=\overline{w}(z)+Cg(z)$ hence
$\underline{v}(z)$ or $\overline{v}(z)$ is determined by solving the
following system
\begin{equation}\label{7}
\begin{array}{ll}
\underline{v}(z)=(B+C)^{-1}(\underline{w}(z)+Cg(z)),\\
\\
\overline{v}(z)=(B+C)^{-1}(\overline{w}(z)+Cg(z)).
\end{array}
\end{equation}
Therefore, the solution of FSLEs $(\ref{1})$ can be obtained by
solving system (\ref{7}) that the vector of solution is unique. But
may still not be an appropriate fuzzy number vector.

\begin{theor}\label{t22}
Let
$\tilde{v}(z)=(\tilde{v}_1(z),\tilde{v}_2(z),\ldots,\tilde{v}_n(z))^{T}$
be the fuzzy solution of Eq. $(\ref{1})$ and the matrix
$\textbf{A}^{-1}$ exists. Then the solution of system
  \begin{equation}\label{3.1}
  \textbf{A}(\overline{\textbf{v}}(z)+\underline{\textbf{v}}(z))=\overline{\textbf{w}}(z)+\underline{\textbf{w}}(z)~,
 \end{equation}
for
$\overline{w}(z)+\underline{w}(z)=(\overline{w}_1(z)+\underline{w}_1(z),\overline{w}_2(z)+\underline{w}_2(z)
,\ldots,\overline{w}_n(z)+\underline{w}_n(z))^{T}$ is in the
following form
$$\overline{v}(z)+\underline{v}(z)=(\overline{v}_1(z)+\underline{v}_1(z),\overline{v}_2(z)+
\underline{v}_2(z),\ldots,\overline{v}_n(z)+\underline{v}_n(z))^{T}.$$

\end{theor}

Since number of operations to solve the $n\times n$ system are less
than the $2n\times 2n$ system, thus Ezzati's method is better in
comparison with Friedman's method. In Theorem 4 of Ezzati's method
\cite{e} the maximum number of multiplication operations (MNMO) were
obtained which had some problems. The following theorem is the
modified version of Theorem 4 in \cite{e}.

\begin{theor}\label{t1}
Assume $n$ is any integer and $n\geq2$.  and denote by $F_{n}$ and
 $E_{n}$ the MNMOs that are required to
calculate
$$v(z)=(\underline{v}_1(z),\underline{v}_2(z),\ldots,\underline{v}_n(z),
-\overline{v}_1(z),-\overline{v}_2(z),\\ \ldots,-\overline{v}_n(z)
)^{T}=S^{-1}w(z)$$ by Friedman's method and
$$v(z)=(\underline{v}_1(z),\underline{v}_2(z),\ldots\underline{v}_n(z),
\overline{v}_1(z),\overline{v}_2(z),\\
\ldots,\overline{v}_n(z))^{T}$$ by Ezzati's method, respectively,
then
\[F_{n}\geq E_{n},\] and $$F_{n}-E_{n}= 2n^{2}.$$
\end{theor}

\textbf{Proof:} Suppose $h_{n}(A)$~is the MNMOs of computing the
matrix $A^{-1}$. Now, we can write
\[\mathbf{S^{-1}}=\left(%
\begin{array}{cc}
  D & E \\
  E & D \\
\end{array}%
\right),\] where
$$
\begin{array}{l}
\displaystyle D=\frac{1}{2}[(B+C)^{-1}+(B-C)^{-1}],\\
\\
\displaystyle E=\frac{1}{2}[(B+C)^{-1}-(B-C)^{-1}].
\end{array}
$$

Therefore, in order to determine $S^{-1}$, computing matrices
$(B+C)^{-1}$~and~$(B-C)^{-1}$ are required. It is clear
that\[h_{n}(S)=h_{n}(B+C)+h_{n}(B-C)=2h_{n}(A).\] Since
$\tilde{v}(z)\in E^{1}$; $\underline{v}(z)$ and $\overline{v}(z)$,
in the simplest case are lines 
 hence
\[F_{n}=2h_{n}(A)+8n^{2}.\]
For computing
$\underline{v}(z)+\overline{v}(z)=(\underline{v}_1(z)+\overline{v}_1(z),\underline{v}_2(z)+\overline{v}_2(z)
,\ldots,\underline{v}_n(z)+\overline{v}_n(z))^{T}$ from Eq.
$(\ref{1})$ and
 $\underline{v}(z)=(\underline{v}_1(z),\underline{v}_2(z),\ldots\underline{v}_n(z))^{T}$~from
Eq. $(\ref{7})$~and according to Ezzati's method, the MNMOs are
$h_{n}(A)+2n^{2}$ and ~$h_{n}(B+C)+4n^{2}$ respectively. Since
$h_{n}(B+C)=h_{n}(A)$ thus
\[E_{n}=2h_{n}(A)+6n^{2},\] and finally $F_{n}-E_{n}=2n^{2}.\square$

\begin{definition}\label{d5}
\cite{e} Assume
$\tilde{v}(z)=\{(\underline{v}_i(z),\overline{v}_i(z)),~1\leq i\leq
n\}$ is the unique solution of Eqs. $(\ref{1}),(\ref{7})$.
 The fuzzy number vector
$\tilde{P}=\{(\underline{p}_i(z),\overline{p}_i(z)),~1\leq i\leq
n\}$ is defined
by\[\underline{p}_i(z)=\min\{\underline{v}_i(z),\overline{v}_i(z),\underline{v}_i(1)\},\]
\[\overline{p}_i(z)=\max\{\underline{v}_i(z),\overline{v}_i(z),\overline{v}_i(1)\},\]
which is called a fuzzy vector solution of Eqs. $(\ref{1})$ and
$(\ref{7})$.
\end{definition}

If $\underline{p}_i(z)=\underline{v}_i(z)$ and
$\overline{p}_i(z)=\overline{v}_i(z),~1\leq i\leq
n,$~then~$\tilde{P}$ is called a strong fuzzy
solution.~Otherwise,~$\tilde{P}$ is called a weak fuzzy solution
which it is not fuzzy linear system's solution and is not always
fuzzy number vector.

\textbf{Remark 2:} We know that Friedman et al.~\cite{k1} and
Ezzati~\cite{e} find two kinds of solutions, which are called the
weak and the strong solutions. The weak solution is not system's
solution and it is not always the fuzzy number vector \cite{a}.
Hence, we do not interest to find weak fuzzy solution. Also, in
these methods the kind of solutions- strong or weak- are determined
only in the end of method and it is one of important faults of these
methods.

In the next section, a novel method
  for solving a $n \times n$ FSLEs is presented. It is observed our method can least the computing error,~because without carrying
   out further computation,~we can determined that the fuzzy linear system,~has no fuzzy number vector solution.

\section{Main Idea}
In this section, a novel and applicable method to solve the FSLEs is
presented. Several theorems and lemmas are illustrated to improve
the ezzati's method \cite{e}. By using these theorems we show the
number of operations of our method are less than the methods of
Ezzati \cite{e} and Friedman \cite{k1}.

\begin{theor}\label{t2}
Suppose the inverse  matrix of \textbf{B+C} exists and
$\tilde{v}(z)=(\tilde{v}_1(z),\tilde{v}_2(z),\ldots,\tilde{v}_n(z))^{T}$
is a fuzzy solution of Eq~$(\ref{1})$.
Then~$\overline{v}(z)-\underline{v}(z)=\left(\overline{v}_1(z)-\underline{v}_1(z),\overline{v}_2(z)-\underline{v}_2(z)
,\ldots,\overline{v}_n(z)-\underline{v}_n(z)\right)^{T}$~  is the
solution of the following system
  \begin{equation}\label{3.1}
    (\textbf{B+C})(\overline{\textbf{v}}(z)-\underline{\textbf{v}}(z))=
    \overline{\textbf{w}}(z)-\underline{\textbf{w}}(z)~,
 \end{equation}
where~$\overline{w}(z)-\underline{w}(z)=(\overline{w}_1(z)-\underline{w}_1(z),\overline{w}_2(z)-\underline{w}_2(z)
,\ldots,\overline{w}_n(z)-\underline{w}_n(z))^{T}.$

\end{theor}

\textbf{Proof:}  Let
$\tilde{v}_j(z)=(\underline{v}_j(z),\overline{v}_j(z)), 1 \leq j
\leq n$ be the parametric form of $\tilde{v}_j$. For positive values
$a_{ij}'$ and $a_{ij}''$ we have
$$
\begin{array}{l}
  a_{ij}=a_{ij}'-a_{ij}'', \\
    \\
  a_{ij}' a_{ij}''=0,
\end{array}
$$
where $a_{ij}, a_{ij}'$ and $a_{ij}''$ are the coefficients of
matrices $A, B$ and $C$ respectively. By presenting the Eq.
$(\ref{1})$ to the parametric form, for~$i=1,2,\ldots,n$ we get
\[(a_{i1}'-a_{i1}'')(\underline{v}_1(z),\overline{v}_1(z))+\ldots+
(a_{in}'-a_{in}'')(\underline{v}_n(z),\overline{v}_n(z))=(\underline{w}_i(z),\overline{w}_i(z)).\]
Hence
\begin{equation}\label{4}
    \begin{array}{cc}
       a_{i1}'\underline{v}_1(z)-a_{i1}''\overline{v}
_1(z)+a_{i2}'\underline{v}_2(z)-a_{i2}''\overline{v}_2(z)+\ldots+a_{in}'\underline{v}_n(z)-a_{in}''\overline{v}_n(z)
=\underline{w}_i(z),
\end{array}
\end{equation}and
\begin{equation}\label{s4}
      \begin{array}{cc}
a_{i1}'\overline{v}_1(z)-a_{i1}''\underline{v}_1(z)+a_{i2}'\overline{v}_2(z)-a_{i2}''
\underline{v}_2(z)+\ldots+a_{in}'\overline{v}_n(z)-a_{in}''\underline{v}_n(z)=\overline{w}_i(z).
\end{array}
\end{equation}
Now, we can differentiate Eq. $(\ref{4})$ from  Eq. $(\ref{s4})$ as
follows
\[(a_{i1}'+a_{i2}'')(\overline{v}_1(z)-\underline{v}_1(z))+(a_{i2}'+a_{i2}'')(\overline{v}_
2(z)-\underline{v}_2(z))+\ldots+(a_{in}'+a_{in}'')(\overline{v}_n(z)-\underline{v}_n(z))
=\overline{w}_i(z)-\underline{w}_i(z).\]
Therefore,~$d(z)=\overline{v}(z)-\underline{v}(z)=(\overline{v}_1(z)-\underline{v}_1(z),\overline{v}_2(z)-\underline{v}_2(z)
,\ldots,\overline{v}_n(z)-\underline{v}_n(z))^{T}$ is the solution
of
$(B+C)(\overline{v}(z)-\underline{v}(z))=\overline{w}(z)-\underline{w}(z)$.$\square$
\begin{theor}\label{t3}
Suppose the inverse matrix of \textbf{B+C} exists.  The
equation~$(\ref{1})$, does not have a fuzzy number vector
solution,~if the vector solution of the following system is not
nonnegative, i.e. at least one of the entries are negative
 \begin{equation}\label{s1}
    (\textbf{B+C})(\overline{\textbf{v}}(z)-\underline{\textbf{v}}(z))=
    \overline{\textbf{w}}(z)-\underline{\textbf{w}}(z).
\end{equation}
\end{theor}
\textbf{Proof:} We know that, the vector solution of the Eq.
$(\ref{s1})$ is ($\overline{v}(z)-\underline{v}(z)$). Now, suppose
that ($\overline{v}(z)-\underline{v}(z)$) is not nonnegative. So,
according to the Definition~$(\ref{d2})$, the fuzzy number vector
solution is not exist. It is clear that the matrix~$(B+C)^{-1}$~is
the non positive matrix, i.e. at least one of the entries is
positive
 because $(B+C)$ is the positive matrix.$\square$


Triangular fuzzy numbers are the simple and  the popular fuzzy
numbers. Also, triangular fuzzy numbers have a special property
$$\overline{w}(z)-\underline{w}(z)=(\rho'+\rho'')(1-z).$$ Hence
when the right hand side vector $\tilde{w}(z)$ is triangular, the
parametric linear system $(\ref{s1})$ can be transformed to the
crisp linear system.

$\mathbf{Lemma~1}$ Suppose the inverse  matrix of $(B+C)$ exists,
and $\tilde{w}(z)\in TE^{1}$. The Eq. $(\ref{1})$,~does not have a
fuzzy number vector solution,~if the vector solution of the
following system is not nonnegative, i.e. at least one of the
entries is negative
\begin{equation}\label{s6}
    (B+C)(\mu'+\mu'')=(\rho'+\rho''),
\end{equation}
where $(\mu'+\mu'')(1-z)=\overline{v}(z)-\underline{v}(z),
~(\rho'+\rho'')(1-z)=\overline{w}(z)-\underline{w}(z).$

\textbf{Proof:} It is clear.$\square$\\

~Now, a new method to solve the FSLEs is presented.~Assume that the
inverse matrix of $A$ in Eq. $(\ref{1})$ exists.~For solving Eq.
$(\ref{1})$, the following system
 \begin{equation}\label{12}
    (\textbf{B+C})(\overline{\textbf{v}}(z)-\underline{\textbf{v}}(z))
    =\overline{\textbf{w}}(z)-\underline{\textbf{w}}(z),
 \end{equation}
should be solved where the matrices $B$ and $C$ were defined in
Subsection 2.2.
    Let the solution of this system be in the following form
$$\mathbf{d(z)}=\left(%
\begin{array}{c}
   d_{1}(z)\\
   \\
   d_{2}(z)\\
   \\
   \vdots\\
   \\
   d_{n}(z)\\
\end{array}%
\right)=
\mathbf{\overline{v}(z)-\underline{v}(z)}=\left(%
\begin{array}{cccc}
   \overline{v}_1(z)-\underline{v}_1(z)\\
   \\
   \overline{v}_2(z)-\underline{v}_2(z)\\
   \\
   \vdots\\
   \\
   \overline{v}_n(z)-\underline{v}_n(z)\\
\end{array}%
\right).
$$

If~$d=\overline{v}-\underline{v}$~is not nonnegative
,~then we do not have the fuzzy number vector
solution. Otherwise, in order to show the existence of fuzzy  number
vector solution for Eq. $(\ref{1})$, we continue our idea. At first,
we should solve the following system
 \begin{equation}\label{12.2}
 \textbf{A}(\overline{\textbf{v}}(z)+\underline{\textbf{v}}(z))=\overline{\textbf{w}}(z)+\underline{\textbf{w}}(z).
 \end{equation}

According to the Theorem $\ref{t22}$, we know that this system has
the solution in the following form
$$\mathbf{g(z)}=\left(%
\begin{array}{c}
   g_{1}(z)\\
   \\
   g_{2}(z)\\
   \vdots\\
   g_{n}(z)\\
\end{array}%
\right)=
\mathbf{\overline{v}(z)+\underline{v}(z)}=\left(%
\begin{array}{cccc}
   \overline{v}_1(z)+\underline{v}_1(z)\\
   \\
   \overline{v}_2(z)+\underline{v}_2(z)\\
   \vdots\\
   \overline{v}_n(z)+\underline{v}_n(z)\\
\end{array}%
\right).
$$
Finally, by solving systems $(\ref{12})$ and $(\ref{12.2})$ and
finding $\mathbf{d(z)}$ and $\mathbf{g(z)}$ we have
\begin{equation}\label{11}
 \left\{\begin{array}{cc}
 \underline{v}(z)=\frac{\mathbf{g(z)-d(z)}}{2},\\
 \\
\overline{v}(z)=\frac{\mathbf{g(z)+d(z)}}{2}.\\
\end{array}
\right.\\
\end{equation}

If the conditions of Definition \ref{d2} are true, then the solution
of FSLEs $(\ref{1})$ can be obtained by solving the crisp linear
system of Eqs. $(\ref{12})$ and $(\ref{12.2})$ that the solution
vector is the fuzzy number vector and unique. Otherwise, if at least
one of the conditions do not true, the fuzzy linear system of Eqs.
(\ref{1}) does not have fuzzy number vector solution.

$\mathbf{Remark~3:}$ If $\tilde{w}\in TE^{1}$, then according to the
Lemma 1 the system of Eqs. (\ref{12}) have the vector solution as
$d'=\mu'+\mu''$ where $d(z)=d'(1-z)$. So, the Eqs. (\ref{11}) can be
written in following form
\begin{equation}\label{11.1}
 \left\{\begin{array}{cc}
 \underline{v}(z)=\frac{\mathbf{g(z)-d'(1-z)}}{2},\\
 \\
\overline{v}(z)=\frac{\mathbf{g(z)+d'(1-z)}}{2}.\\
\end{array}
\right.\\
\end{equation}

\begin{theor}\label{t4}
  Assume that $n$ is any integer, $n\geq2$ and denote by
  $E_{n}$~and~$D_{n}$
the MNMOs that are required to calculate
 $$v(z)=(\underline{v}_1(z),\underline{v}_2(z),\ldots\underline{v}_n(z),\overline{v}_1(z),
 \overline{v}_2(z),\ldots,\overline{v}_n(z))^{T},$$
 in the Ezzati's method~$\cite{e}$ and presented method then

  \begin{equation}\label{scalar}
 \left\{\begin{array}{lc}
 E_{n}-D_{n}=2n^{2}, & \overline{v}(z)-\underline{v}(z)\geq0,\\
 \\
 E_{n}-D_{n}=h_{n}(A)+4n^{2},& o.w.,\\
\end{array}
\right.
\end{equation}
where $h_{n}(A)$ shows the MNMOs that are required to calculate
$A^{-1}$.
  \end{theor}
$\mathbf{Proof:}$
  According to the Theorem~$\ref{t1}$,~we have
  \[E_{n}=2h_{n}(A)+6n^{2}.\]
Assume $d(z)=\overline{v}(z)-\underline{v}(z)$~is the nonnegative
  matrix 
  then
  for understanding that whether the fuzzy linear system of Eqs. $(\ref{1})$~has the fuzzy
  number vector solution, we need to solve the system of Eqs. $(\ref{12.2})$.~So,~for
  computing $$\overline{v}(z)-\underline{v}(z)=(\overline{v}_1(z)-\underline{v}_1(z),\overline{v}_2(z)-\underline{v}_2(z)
,\ldots,\overline{v}_n(z)-\underline{v}_n(z))^{T},$$ and
 $$\overline{v}(z)+\underline{v}(z)=(\overline{v}_1(z)+\underline{v}_1(z),\overline{v}_2(z)+\underline{v}_2(z)
,\ldots,\overline{v}_n(z)+\underline{v}_n(z))^{T},$$ from Eqs.
$(\ref{12})$ and $(\ref{12.2})$, the maximum number of
multiplication operations are $h_{n}(B+C)+2n^{2}$~and
 $h_{n}(A)+2n^{2}$,~respectively.~Clearly $h_{n}(B+C)=h_{n}(A)$.
 Hence
\[D_{n}= 2h_{n}(A)+4n^{2},\]
and
\[E_{n}-D_{n}=2n^{2}.\]

Otherwise, assume that~$d(z)=\overline{v}(z)-\underline{v}(z)$~is
not the nonnegative matrix. According to the Theorem $\ref{t3}$ we
do not have the fuzzy number vector solution for solving the FSLEs
$(\ref{1})$. If we do not have the fuzzy number vector
solution,~there will be no necessary for computing
~$\overline{v}(z)+\underline{v}(z)$~from system of Eqs.
$(\ref{12.2})$. Thus,~we need to compute
$d(z)=\overline{v}(z)-\underline{v}(z)$. Therefore,~in this case,
the maximum number of multiplication operations are
$h_{n}(B+c)+2n^{2}$~or~$h_{n}(A)+2n^{2}$. Finally, we have

\[E_{n}\geq D_{n},~~~E_{n}-D_{n}=h_{n}(A)+4n^{2}. \square \]
$\mathbf{Lemma~2}$ ~Let~$\tilde{w}(z)$ be the triangular fuzzy
number vector,~from Eq. $(\ref{1})$. Then~$\tilde{v}(z)$ is the
triangular fuzzy number vector solution,~from Eq. $(\ref{1})$.

$\mathbf{Proof:}$~It is clear. $\square$\\
$\mathbf{Lemma~3}$ Suppose that in the
Theorem~$\ref{t4},~\tilde{w}(z)$~is the triangular fuzzy number
vector from Eq. $(\ref{1})$,~then~$E_{n}\geq D_{n}$~and
  \begin{equation}\label{scalar}
 \left\{\begin{array}{cc}
 E_{n}-D_{n}=3n^{2}-n,~~~~~&\overline{v}(z)-\underline{v}(z)\geq0,\\
 \\
 E_{n}-D_{n}=h_{n}(A)+5n^{2},&o.w.\\
\end{array}
\right.\\
\end{equation}\\

$\mathbf{Proof:}$ If~$\tilde{w}(z)$~is the triangular fuzzy number
vector from Eq. $(\ref{1})$~i.e., since $\tilde{w}(z)\in TE^{1};
\underline{w}(z)$ and $\overline{w}(z)$ in the simplest case is the
line. So clearly,~according to the Theorem~$\ref{t1}$,~we have
\[E_{n}=2h_{n}(A)+6n^{2},\]
and according to the Remark 1,~we get
\[\overline{v}(z)=c+(1-z)\mu',~~~~~\underline{v}(z)=c-(1-z)\mu''.\]
So,
\[\overline{v}(z)-\underline{v}(z)=(\mu'+\mu'')(1-z),\]
 ~\[\overline{w}(z)-\underline{w}(z)=(\rho'+\rho'')(1-z),\]

 and from the system of Eqs. $(\ref{12})$,~we have
 \[(B+C)(\mu'+\mu'')(1-z)=(\rho'+\rho'')(1-z).\]
 If ~$r\neq1$, the following relation can be obtained as
\begin{equation}\label{s5}
    (B+C)(\mu'+\mu'')=\rho'+\rho'',
\end{equation} where it is the crisp linear system. It is clear that for $r=1$ the FSLE replaced by crisp linear system.

Now,~we assume that~$d'=(\mu'+\mu'')$~is nonnegative. Then for
understanding that whether the FSLEs $(\ref{1})$~has the fuzzy
number vector solution, we need to solve the system of Eqs.
$(\ref{12.2})$. So,~for computing
$(\mu'+\mu'')=(\mu'_1+\mu''_1,\mu'_2+\mu''_2,\ldots,\mu'_n+\mu''_n)$
~from Eq. (\ref{s5}) and
$\overline{v}(z)+\underline{v}(z)=(\overline{v}_1(z)+\underline{v}_1(z),\overline{v}_2(z)+\underline{v}_2(z)
,\ldots,\overline{v}_n(z)+\underline{v}_n(z))^{T}$~from Eq.
(\ref{12.2}) and $d(z)=(1-z)d'$ for final solution in Eq.
(\ref{11.1})
  the maximum number of multiplication operations are $h_{n}(B+C)+n^{2}$, $h_{n}(A)+2n^{2}$ and $n$
  respectively. Clearly $h_{n}(B+C)=h_{n}(A)$. So \[D_{n}=2h_{n}(A)+3n^{2}+n,\]
  and \[E_{n}\geq D_{n},~~~E_{n}-D_{n}=3n^{2}-n.\]
Otherwise,~assume that~$d~'=(\mu'+\mu'')$~ is not nonnegative. Then,
according to the Lemma~1 we do not have a fuzzy number vector
solution for solving the fuzzy linear system of Eqs. $(\ref{1})$. We
know that if we do not have a fuzzy number vector solution,~there
will be no necessary for computing
~$\underline{v}(z)+\overline{v}(z)$~from Eq. $(\ref{12.2})$.
Thus,~we need to compute~$d~'=(\mu'+\mu'')$.
 Therefore,~in this case the maximum number of multiplication operations are $$h_{n}(B+C)+n^{2}~or~h_{n}(A)+n^{2}.$$ Then,
 ~we have\[E_{n}\geq D_{n},~~~E_{n}-D_{n}=h_{n}(A)+5n^{2}.\square\]


\section{Numerical illustrations}
In this section, some examples of the FSLEs are presented \cite{k1}.
Also, two algorithms are applied to solve the problems. Furthermore,
several graphs are demonstrated that show the fuzzy form of
solutions.


\textbf{Algorithm 1:} Let $\textbf{A}$ be the nonsingular matrix.\\
\texttt{Step 1: Input matrix~$A=[a_{ij}]\in R^{n\times n}$
and~$\tilde{v}(z)=(\underline{v}(z),\overline{v}(z)),
\tilde{w}(z)=(\underline{w}(z),\overline{w}(z))\in E^{1}$.}\\
\texttt{Step 2: Calculate~$B=[b_{ij}]$ and $C=[c_{ij}]$.\\
$$
\left\{
\begin{array}{l}
If~a_{ij}>0~\Rightarrow~b_{ij}=a_{ij};
                ~else~b_{ij}=0, \\
   \\
If~a_{ij}<0~\Rightarrow~ c_{ij}=-a_{ij};~
                 else~ c_{ij}=0.
\end{array}\right.
$$}\\
\texttt{Step 3: Calculate $M=(B+C)^{-1}$.}\\
\texttt{Step 4: Calculate
$d_{1}(z)=\overline{w}(z)-\underline{w}(z)$ and
 $d(z)=M.d_{1}(z)$.  If  $d(z)$  is not nonnegative,}

\texttt{ go to step 8.}\\
\texttt{Step 5: Calculate $k=A^{-1},
g_{1}(z)=\overline{w}(z)+\underline{w}(z)$ and
                 $g(z)=k.g_{1}(z)$.}\\
\texttt{Step 6: Calculate $\underline{v}(z)=\frac{g(z)+d(z)}{2}$ and
$\overline{v}(z)=\frac{g(z)-d(z)}{2}$. }\\
\texttt{Step 7: If conditions of Definition \ref{d2} are true then
$\tilde{v}(z)=(\underline{v}(z),\overline{v}(z))$ and go to}

\texttt{step 9. Else go to step 8.}\\
\texttt{Step 8: Show the message "The system does not have fuzzy number vector solution".}\\
 \texttt{Step 9: End.}\\

The following algorithm is presented to triangular fuzzy linear
system.

\textbf{Algorithm 2:}\\
\texttt{Step 1: Input matrix~$A=[a_{ij}]\in R^{n\times n}$
and~$\tilde{v}(z)=(\underline{v}(z),\overline{v}(z)),
\tilde{w}(z)=(\underline{w}(z),\overline{w}(z))\in TE^{1}$.}\\
\texttt{Step 2: Calculate~$B=[b_{ij}]$ and $C=[c_{ij}]$.\\
$$ \left\{
\begin{array}{l}
If~a_{ij}>0~\Rightarrow~b_{ij}=a_{ij};
                ~else~b_{ij}=0, \\
   \\
If~a_{ij}<0~\Rightarrow~ c_{ij}=-a_{ij};~
                 else~ c_{ij}=0.
\end{array}\right.
$$}\\
\texttt{Step 3: Calculate $M=(B+C)^{-1}$.}\\
\texttt{Step 4: Calculate~$d_{1}=\rho'+\rho''$ and $d'=M.d_{1}$. If
$d'$ is not nonnegative,~go to}

\texttt{step 8.}\\
\texttt{Step 5: Calculate~$k=A^{-1}, g_{1}(z)=\overline{w}(z)+\underline{w}(z)$ and $g(z)=k.g_{1}(z)$.}\\
\texttt{Step 6: Calculate~$\underline{v}(z)=\frac{g(z)+d'(1-z)}{2};$
and $\overline{v}(z)=\frac{g(z)-d'(1-z)}{2}$.}\\
\texttt{Step 7: If conditions of Definition \ref{d2} are true then
$\tilde{v}(z)=(\underline{v}(z),\overline{v}(z))$ and go to}

\texttt{step 9. Else go to step 8.}\\
\texttt{Step 8: Show the message "The system does not have fuzzy
number vector
solution".}\\
\texttt{  Step 9: End.}


\begin{example}  \cite{k1} Consider the following $2\times 2$ FSLEs
$$
\left\{\begin{array}{l}
\tilde{v}_1-\tilde{v}_2=(z,2-z),\\
\\
\tilde{v}_1+3\tilde{v}_2=(4+z,7-2z),\\
\end{array}
\right.
$$
where ~$\tilde{w}$ is a triangular vector of fuzzy numbers hence the
Algorithm 2 is applied. By using this algorithm we have:\\
\textbf{Step1.} Input matrix $$  A=\left(%
                                                            \begin{array}{cc}
                                                                 1&-1 \\
                                                                 \\
                                                                1&3 \\
                                                                \end{array}%
                                                                    \right),$$
                                                                    $$\tilde{v}(z)=
                                                                    (\underline{v}(z),\overline{v}(z))
                                                                    \in
                                                                    TE^{1},$$
                                                                    $$\tilde{w}(z)=
                                                                    (\underline{w}(z),\overline{w}(z))=\left(%
                                                            \begin{array}{c}
                                                             (z,2-z) \\
                                                             \\
                                                             ( 4+z,7-2z)\\
                                                            \end{array}
                                                            \right).
                                                            $$\\
                \textbf{Step 2.} Calculate $$B=[b_{ij}]=\left(%
                        \begin{array}{cc}
                            1&0 \\
                            \\
                            1&3\\
                            \end{array}
                            \right),$$ $$C=[c_{ij}]=\left(%
                        \begin{array}{cc}
                            0&1 \\
                            &\\
                            0&0\\
                            \end{array}
                            \right).$$\\
                 \textbf{Step~3.} Compute $$M=(B+C)^{-1}=\left(%
                        \begin{array}{cc}
                         \displaystyle   \frac{3}{2}&  \displaystyle-\frac{1}{2} \\
                           & \\
                \displaystyle   -\frac{1}{2}& \displaystyle \frac{1}{2}\\
                            \end{array}
                            \right).$$\\
                \textbf{Step~4.} Calculate $$d_{1}=\rho'+\rho''=\left(%
                                                        \begin{array}{c}
                                                         \rho'_1+\rho''_1 \\
                                                         \\
                                                          \rho'_2+\rho''_2 \\
                                                        \end{array}%
                                                        \right)=\left(%
                                                        \begin{array}{cc}
                                                           2\\
                                                           \\
                                                          3\\
                                                        \end{array}%
                                                        \right),$$ $$d'=M.d_{1}=\left(%
                                                    \begin{array}{c}
                                            \displaystyle   \frac{3}{2}\\
                                                        \\
                                              \displaystyle \frac{1}{2}\\
                                                        \end{array}
                                                        \right).$$

                                                        It is clear that~$d'=\mu'+\mu''$~is the nonnegative matrix,~therefore
                           go to the step 5.\\
                \textbf{Step 5.} Calculate
                $$~k=A^{-1}=\left(%
                                                    \begin{array}{cc}
                                               \displaystyle   \frac{3}{4}& \displaystyle \frac{1}{4}\\
                                               \\
                                             \displaystyle  -\frac{1}{4}& \displaystyle \frac{1}{4}\\
                                                        \end{array}%
                                                        \right),$$
                                                       $$g_{1}(z)=\overline{w}(z)+\underline{w}(z)=\left(%
                                                    \begin{array}{c}
                                                        2\\
                                                        \\
                                                        11-z\\
                                                        \end{array}%
                                                        \right),$$
                                                        $$g(z)=k.g_{1}(z)=\left(%
                                                    \begin{array}{cccccc}
                                                         \displaystyle \frac{17-z}{4}\\
                                                         \\
                                                         \displaystyle \frac{9-z}{4}\\
                                                        \end{array}%
                                                        \right).$$\\
                                                   \textbf{Step~6.} Compute   $$
                 \underline{v}(z)=\frac{g(z)-d(z)}{2}=\left(%
                                                    \begin{array}{c}
                                                      \underline{v}_1 \\
                                                      \\
                                                      \underline{v}_2 \\
                                                    \end{array}%
                                                    \right)=\left(%
                                                    \begin{array}{c}
                                                      1.375+0.625z\\
                                                      \\
                                                      0.875+0.125z \\
                                                    \end{array}%
                                                    \right),
                                                    \\$$
                                                     $$ \overline{v}(z)=\frac{g(z)+d(z)}{2}=\left(%
                                                        \begin{array}{cc}
                                                          \overline{v}_1 \\
                                                          \\
                                                          \overline{v}_2 \\
                                                        \end{array}%
                                                        \right)=\left(%
                                                        \begin{array}{cccc}
                                                          2.875-0.875z\\
                                                          \\
                                                          1.375-0.375z \\
                                                        \end{array}%
                                                        \right).\\$$

\begin{figure}
\centering
$$\begin{array}{c}
  \includegraphics[width=3in]{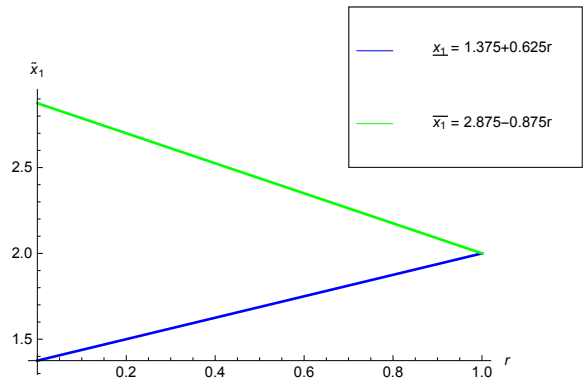} \\
    \\
    \\
    \\
  \includegraphics[width=3in]{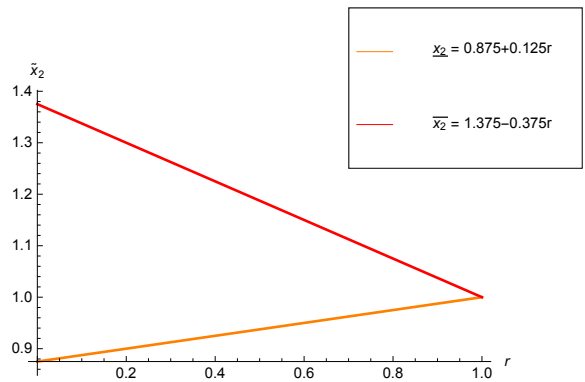}
\end{array}$$
  \\
  \caption{The solutions of Example 1.}\label{F1}
\end{figure}

                                                     Since the conditions of Definition
                                                     \ref{d2}
                                                      are true, the FSLEs has the fuzzy number vector
                                                      solution. Fig. \ref{F1}
                                                      shows the
                                                      obtained solution is in the fuzzy
                                                      form.
                                                      \end{example}

\begin{figure}[h]
\centering
$$\begin{array}{c}
  \includegraphics[width=3in]{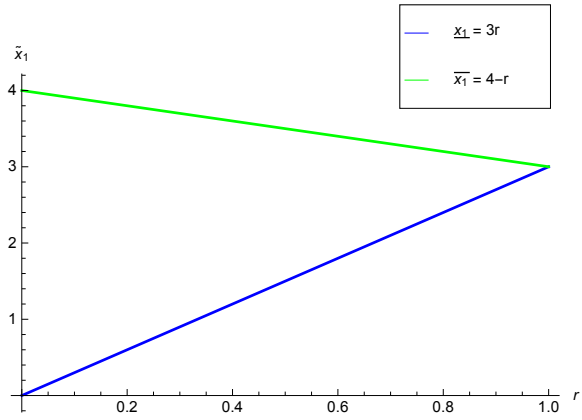} \\
    \\
    \\
    \\
  \includegraphics[width=3in]{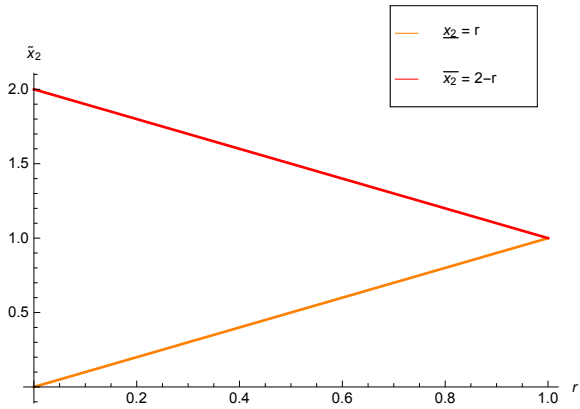}
\end{array}$$
  \\
  \caption{The solutions of Example 2.}\label{F2}
\end{figure}

\begin{example}~\cite{k1}~
 Consider the~$3\times 3$ FSLEs
 \[
 \left\{\begin{array}{l}
\tilde{v}_1+\tilde{v}_2-\tilde{v}_3=(z,2-z),\\
\\
\tilde{v}_1-2\tilde{v}_2+\tilde{v}_3=(2+z,3),\\
\\
2\tilde{v}_1+\tilde{v}_2+3\tilde{v}_3=(-2,-1-z), \\
\end{array}
\right.
\]
where $\tilde{Y}$~is a triangular fuzzy number vector. By using
Algorithm~2 we  have\\
\textbf{Step~1.} Input matrix $$            A=\left(%
                                                            \begin{array}{ccc}
                                                                 1&1&-1 \\
                                                                 &&\\
                                                                1&-2&~~1 \\
                                                                &&\\
                                                            2& 1& 3\\
                                                                \end{array}%
                                                                    \right),$$

                                                                    $$
                                 \tilde{w}(z)=(\underline{w}(z),\overline{w}(z))=\left(%
                                                            \begin{array}{c}
                                                            (z,2-z) \\
                                                            \\
                                                            (2+z,3)\\
                                                            \\
                                                             (-2,-1-z)\\
                                                            \end{array}
                                                            \right).$$\\
               \textbf{Step~2.} Compute $$B=[b_{ij}]=\left(%
                        \begin{array}{ccc}
                          1&1&0 \\
                          &&\\
                             1&0&1 \\
                             &&\\
                          2&1&3\\
                            \end{array}
                            \right),$$
                            $$C=[c_{ij}]=\left(%
                        \begin{array}{ccc}
                            0&0&1 \\
                           && \\
                             0&2&0 \\
                             &&\\
                          0&0&0\\
                            \end{array}
                            \right).$$\\
                \textbf{Step~3.} Calculate  $$M=(B+C)^{-1}=\left(%
                        \begin{array}{ccc}
                            5&-1&-3 \\
                            \\
                             -2& 1& 1 \\
                             \\
                          -1& 0& 1\\
                            \end{array}
                            \right).$$\\
                \textbf{Step~4.} Compute $$d_{1}=\rho'+\rho''=\left(%
                                                        \begin{array}{c}
                                                         \rho'_1+\rho''_1 \\
                                                         \\
                                                          \rho'_2+\rho''_2 \\
                                                        \end{array}%
                                                        \right)=\left(%
                                                        \begin{array}{c}
                                                          2\\
                                                          \\
                                                          1\\
                                                          \\
                                                          1\\
                                                        \end{array}%
                                                        \right),$$
                                             $$ d'=M.d_{1}=\left(%
                                                    \begin{array}{c}
                                                          6\\
                                                          \\
                                                          -2 \\
                                                          \\
                                                          -1\\
                                                        \end{array}
                                                        \right).$$
                                                        Since $d'=\mu'+\mu''$ is not nonnegative, therefore
                           the system has not a fuzzy number vector solution.

                                                      \end{example}
\begin{example}
Consider the following $2\times 2$ FSLEs
\begin{equation}\label{scalar}
\left\{\begin{array}{cc}
\tilde{v}_1+\tilde{v}_2=(4z,6-2z),\\
\\
\tilde{v}_1+2\tilde{v}_2=(5z,8-3z).\\
\end{array}
\right.\end{equation} In this example, by applying the Algorithm 1
we have\\
\textbf{Step 1.} Input matrix $$A=\left(%
\begin{array}{cc}
1&1 \\
&\\
1&2 \\
\end{array}%
\right),$$  $$\tilde{w}(z)=(\underline{w}(z),\overline{w}(z))=\left(%
\begin{array}{c}
(4z,6-2z) \\
\\
(5z,8-3z)\\
\end{array}
\right).$$\\
\textbf{Step~2.} Compute $$B=[b_{ij}]=\left(%
 \begin{array}{cc}
 1&1\\
 &\\
 1&2\\
 \end{array}
 \right),$$ $$C=[c_{ij}]=\left(%
 \begin{array}{cc}
  0&0\\
  &\\
  0&0\\
  \end{array}
  \right).$$\\
\textbf{Step~3.} Calculate $$M=(B+C)^{-1}=\left(%
\begin{array}{cc}
2&-1 \\
\\
-1&1\\
 \end{array}
 \right).$$\\
\textbf{Step~4.} Compute $$d_{1}(z)=\overline{w}(z)-\underline{w}(z)=\left(%
 \begin{array}{c}
 6-6z \\
 \\
 8-8z\\
  \end{array}%
  \right),$$ $$d(z)=M.d_{1}=\left(%
    \begin{array}{c}
    4-4z\\
    \\
    2-2z\\
     \end{array}
     \right).$$

Since ~$d(z)=\overline{v}(z)-\underline{v}(z)$~is nonnegative for
~$0\leq z\leq1$ ,therefore, in order to find that whether the fuzzy
system has the fuzzy number vector solution we need to go to step
5.\\
 \textbf{Step~5.}  Calculate
                $$k=A^{-1}=\left(%
                                                    \begin{array}{cc}
                                                        2&-1\\
                                                        &\\
                                                        -1&1\\
                                                        \end{array}%
                                                        \right),$$
                                                       $$g_{1}(z)=\overline{w}(z)+\underline{w}(z)=\left(%
                                                    \begin{array}{c}
                                                        6+2z\\
                                                        \\
                                                        8+2z\\
                                                        \end{array}%
                                                        \right),$$
                                                        $$g(z)=k.g_{1}(z)=\left(%
                                                    \begin{array}{c}
                                                        4+2z\\
                                                        \\
                                                        2\\
                                                        \end{array}%
                                                        \right).$$\\
                                                       \textbf{Step~6.
                                                       } Compute
                                                        $$\underline{v}(z)=\frac{g(z)-d(z)}{2}=\left(%
 \begin{array}{c}
 3z\\
 \\
 z\\
  \end{array}%
  \right),$$
     $$ \overline{v}(z)=\frac{g(z)+d(z)}{2}=\left(%
  \begin{array}{c}
  4-z \\
  \\
  2-z \\
   \end{array}%
 \right).$$

Since the conditions of Definition \ref{d2} are connected hence the
vector solution is the fuzzy number vector solution. Thus,
   the FSLEs (\ref{scalar}) has the fuzzy number vector
   solution. Fig. \ref{F2} shows the fuzzy solutions of this
   example.
\end{example}

\section{Conclusion}
Weakly fuzzy solution was introduced by Friedman et al. \cite{k1}.
This solution is not always fuzzy number vector
 and it is not fuzzy linear system's solution. Also, the kind of solution is only determined in the
 end of solving problem. Hence, it is important to introduce a novel method for solving the FSLEs and find its fuzzy number vector solution.
In the novel proposed method, the original fuzzy system is replaced
by two $n\times n $ crisp linear system. By proving several theorems
we showed the number of operations for presented method are less
than the methods of Friedman and Ezzati. Presented algorithms show
the accuracy and efficiency of method to solve the examples.


\section*{Acknowledgements}
The work of J. J. Nieto has been partially supported by Agencia
Estatal de Investigaci\'on (AEI) of Spain under grant
MTM2016-75140-P, co-financed by the European Community fund FEDER,
and XUNTA de Galicia under grants GRC2015-004 and R2016-022. 



\end{document}